\documentclass[11pt]{amsart}
\usepackage{amssymb}

\newtheorem{theorem}{Theorem}[section]
\newtheorem{prop}[theorem]{Proposition}
\newtheorem{proposition}[theorem]{Proposition}
\newtheorem{lemma}[theorem]{Lemma}

\newtheorem{corollary}[theorem]{Corollary}
\theoremstyle{definition}
\newtheorem{definition}[theorem]{Definition}
\newtheorem{example}[theorem]{Example}

\theoremstyle{remark}
\newtheorem{remark}[theorem]{Remark}
\numberwithin{equation}{section}
\def\sym{{\rm Sym}}
\def\Hom{{\rm Hom}}
\def\counit{\varepsilon}
\def\sign{\epsilon}

\author[Thomas Lam]{Thomas Lam}
\address{Department of Mathematics, Harvard University, Cambridge MA, 02138, USA}
\email{tfylam@math.harvard.edu}

\urladdr{http://www.math.harvard.edu/~tfylam}

\title{Schubert polynomials for the affine Grassmannian}

\date{March 2006}
\thanks{I am indebted to my coauthors Luc Lapointe,
Jennifer Morse and Mark Shimozono, with whom I have studied
$k$-Schur functions.  I began working on $k$-Schur and dual
$k$-Schur functions more than a year ago when Jennifer first
introduced them to me, and Mark explained his geometric conjectures
to me.  I would also like to thank Shrawan Kumar for comments on an
earlier version of this paper.}

\def\ll{\lambda}
\def\N{\mathbb N}
\def\Q{\mathbb Q}
\def\Frac{{\rm Frac}}
\def\O{{\mathbf{O}}}
\def\F{{\mathbf{F}}}
\def\Gr{{\mathcal Gr}}
\newcommand\cd{{\rm cd}}

\newcommand{\ip}[1]{\left<#1\right>}
\newcommand\Z{\mathbb Z}

\newcommand\C{\mathbb C}
\newcommand\A{\mathbb A}
\newcommand\B{\mathbb B}
\newcommand\BB{{\mathcal B}}
\newcommand\G{{\mathcal G}}
\newcommand\p{{\mathcal P}}

\newcommand\id{\mathrm{id}}
\newcommand\aff{\mathrm{aff}}
\newcommand\tF{\tilde{F}}

\begin{document}
\begin{abstract} Confirming a conjecture of Mark Shimozono, we
identify polynomial representatives for the Schubert classes of the
affine Grassmannian as the $k$-Schur functions in homology and
affine Schur functions in cohomology.  Our results rely on the study
of certain subalgebras of Kostant and Kumar's nilHecke ring, and
certain distinguished elements called {\it non-commutative $k$-Schur
functions}.  We also use work of Peterson on the homology of based
loops on a compact group.
\end{abstract}
\maketitle
\section{Introduction}
Let $G$ be a complex simply connected simple group and $K$ a maximal
compact subgroup.  Let $\Gr = G(\C((t)))/G(\C[[t]])$ denote the
affine Grassmannian. The space $\Gr$ is homotopy equivalent to the
space of based loops $\Omega K$ into the compact group $K$.  This
gives the homology $H_*(\Gr)$ and cohomology $H^*(\Gr)$ the
structures of dual Hopf-algebras.

The Bruhat decomposition of $G(\C((t)))$ induces a stratification of
$\Gr$ by Schubert cells, indexed by {\it Grassmannian} elements $w$
of the affine Weyl group $W_\aff$.  Our main theorem
(Theorem~\ref{thm:main}) is the explicit identification, as
polynomials, of the Schubert classes $\sigma_{w} \in H_*(\Gr)$ and
$\sigma^w \in H^*(\Gr)$ when $G = SL(n,\C)$. These polynomials are
combinatorially defined symmetric functions: in homology they are
Lascoux-Lapointe-Morse's $k$-Schur functions $s^{(k)}_w(x)$ and in
cohomology the dual $k$-Schur functions $\tF_w(x)$ (or affine Schur
functions); see~\cite{LLM,LM05,Lam}. Our theorem was originally
conjectured by Mark Shimozono (the conjecture was made explicit in
the cohomology case by Jennifer Morse).

\medskip

The Hopf algebras $H_*(\Gr)$ and $H^*(\Gr)$ were first identified by
Bott~\cite{Bot} using the structure of $\Gr$ as the based loops on a
compact group (another description is given by Ginzburg~\cite{Gin},
which we shall not use).  Separately, Kostant-Kumar~\cite{KK}
studied the topology of homogeneous spaces for arbitrary Kac-Moody
groups, and in particular calculated the structure constants of
$H^*(\Gr)$ in the Schubert basis, using the algebraic construction
of the {\it nilHecke ring}. Our connection between the topology and
the combinatorics proceeds via the study of three subalgebras of the
nilHecke ring.  The first subalgebra is the {\it nilCoxeter algebra}
$\A_0$, which is the algebra generated by the ``divided difference''
operators.  The second algebra is a certain centralizer algebra
$Z_{\A_\aff}(S)$ which we call the {\it Peterson subalgebra}.
Peterson~\cite{Pet} constructs an isomorphism $j$ from the
$T$-equivariant homology $H_T(\Gr)$ of the affine Grassmannian to
$Z_{\A_\aff}(S)$.  The third algebra $\B$ is a combinatorially
defined subalgebra of $\A_0$, which we call the {\it (affine)
Fomin-Stanley subalgebra}. It is (in the finite case) implicitly
used in the construction of Schubert polynomials and Stanley
symmetric functions by Fomin and Stanley~\cite{FS}.

These algebras are tied together by the study of elements $s_w^{(k)}
\in \B$, introduced in~\cite{Lam}, called {\it non-commutative
$k$-Schur functions}.  We show that they have a trio of
descriptions:
\begin{itemize}
\item
they are images of the $k$-Schur functions under an explicit
isomorphism between symmetric functions and the Fomin-Stanley
subalgebra;
\item
they are characterized by a leading ``Grassmannian term'' and a
certain commutation property with the scalars $S = H^T({\rm pt})$;
and
\item
they are the evaluation at 0 of the images of the equivariant
homology Schubert classes $\sigma_{(w)} \in H_T(\Gr)$ under
Peterson's map $j$.
\end{itemize}
The three descriptions establishes the connection between the
topology, algebra and combinatorics.

\medskip

As a corollary of our main result, we establish conjectured
positivity properties of $k$-Schur functions and affine Schur
functions, via geometric positivity results of Graham, Kumar and
Peterson~\cite{Gra, Kum, Pet}.  In the other direction, certain
conjectured combinatorial properties of $k$-Schur functions suggest
interesting geometric properties, the most striking of which is the
connection with Macdonald polynomials.  Our results also gives new
understanding to properties of symmetric functions. For example, the
Hall inner product is given an interpretation as a pairing between
homology and cohomology; Stanley symmetric functions are given a
direct geometric interpretation (Remark~\ref{rem:stan}) and their
symmetry is explained by the commutativity of $H_*(\Gr)$.  We
recover the Schur-positivity of Stanley symmetric
functions~\cite{EG,LS85} and prove the more general (affine
Schur-)positivity for affine Stanley symmetric functions conjectured
in~\cite{Lam}.

Together work with Lapointe, Morse and Shimozono~\cite{LLMS}, we
develop a combinatorial framework for $k$-Schur functions and affine
Schur functions suitable for the study of the affine Grassmannian.
In particular, we obtain an affine Pieri rule.

Some of our results generalize to all simple Lie types, with the
exception of the combinatorics.  We hope to rectify this in a later
work.

\bigskip

{\bf Organisation.}  In Sections~\ref{sec:equi} - \ref{sec:FShom},
we will work in arbitrary Lie type.  In Section~\ref{sec:equi}, we
establish basic notation and facts about the affine Weyl group and
equivariant (co)homology of the affine Grassmannian; we work only in
the case of a simple and simply-connected group.  In
Section~\ref{sec:nilhecke}, we introduce Kostant-Kumar's nilHecke
ring and state the main theorem connecting it with the topology.  In
Section~\ref{sec:pet}, we give a short but hopefully sufficiently
detailed exposition of Peterson's $j$-homomorphism. Since Peterson's
work~\cite{Pet} has yet to be published, we give more explanations
here than just those needed to state Peterson's result. However, we
do not aim to give a systematic development of Peterson's
construction.  In Section~\ref{sec:FShom}, we study a subalgebra of
the nilHecke ring which we call the affine Fomin-Stanley subalgebra.
We show that it is a model for the homology $H_*(\Gr)$.

From Section~\ref{sec:symkschur}, onwards, we specialize to type
$A$. In Section~\ref{sec:symkschur}, we establish notation for
symmetric functions, and following~\cite{Lam}, define affine Stanley
symmetric functions, $k$-Schur functions and non-commutative
$k$-Schur functions.  In Section~\ref{sec:main}, we state and prove
our main theorem. The methods used here are essentially
combinatorial. In Section~\ref{sec:pos}, we explain a number of
positivity properties. In Section~\ref{sec:final}, we give
directions for further study and also identify and explain the
center of the affine nilHecke ring.

\medskip
A preliminary shortened version of this work appeared
as~\cite{LamAS}.  We should caution that some of the notation has
changed.  In particular $W$ now refers to the finite Weyl group and
$W_\aff$ refers to the affine Weyl group.

\section{Equivariant homology and cohomology of the affine Grassmannian}
\label{sec:equi} More details for the material in this section can
be found in~\cite{Kum, Hum}.

\subsection{Affine Weyl group}
Let $W$ be a crystallographic Coxeter group and let $\{r_i \mid i
\in I\}$ denote its simple generators.  Let $R$ be the root system
for $W$.  Let $R^+$ denote the positive roots and $\{ \alpha_i \mid
i \in I\}$ denote the simple roots. Let $Q = \oplus_{i \in I}
\alpha_i$ denote the root lattice and let $Q^\vee = \oplus_{i \in I}
\alpha_i^\vee$ denote the co-root lattice.  Let $h_\Z^*$ and $h_\Z$
denote a weight lattice and co-weight lattice respectively.  We
assume that maps $Q \to h_\Z^*$ and $Q^\vee \to h_\Z$ are given and
fixed.  Let $\ip{.,.}$ denote the pairing between $h_\Z$ and
$h_\Z^*$.

Let $W_\aff = W \ltimes Q^\vee$ denote the affine Weyl group.  Let
$r_0$ denote the additional simple generator of $W_\aff$.  For an
element $\lambda \in Q^\vee$, we let $t_\lambda$ denote the
corresponding translation element in $W_\aff$.  Note that
translations are written multiplicatively: $t_\lambda \cdot t_\mu$ =
$t_{\lambda+\mu}$.  We have the conjugation formula $w\, t_\lambda
w^{-1} = t_{w \cdot \lambda}$. For a real root $\alpha$ of the
affine root system $Q_\aff = \oplus_{i \in I \cup \{0\}} \alpha_i$,
we let $r_\alpha$ denote the corresponding reflection.

Let $\ell: W_\aff \rightarrow \N$ denote the length function of
$W_\aff$ and $\sign(w) = (-1)^{\ell(w)}$ denote the sign.  We note
the following formula for the length of the element $w\, t_\lambda
\in W_\aff$ (where $w \in W$ and $\lambda \in Q^\vee$):
\begin{equation}
\label{eq:length} \ell(w\,t_\lambda) = \sum_{\alpha \in R^+} |
\ip{\lambda,\alpha} + \chi(w\cdot\alpha)|,
\end{equation}
where $\chi(\alpha) = 0$ if $\alpha \in R^+$ and $\chi(\alpha) = 1$
otherwise.

Let $W^0$ denote the minimal length coset representatives of
$W_\aff/W$, which we call {\it Grassmannian} elements. There is a
natural bijection between $W^0$ and $Q^\vee$: each coset $W_\aff/W$
contains one element from each set.  Using the length formula one
can show that $W^0 \cap Q^\vee = Q^-$, the (translations
corresponding to) elements of the co-root lattice which are
anti-dominant.  In fact an element $w t_\lambda$ lies in $W^0$ if
and only if $t_\lambda \in Q^-$ and $w \in W^\lambda$ where
$W^\lambda$ is the set of minimal length representatives of
$W/W_\lambda$ and $W_\lambda$ is the stabilizer subgroup of
$\lambda$.

\subsection{Affine Grassmannian}
Let $G$ be a simple and simply-connected complex algebraic group
with Weyl group $W$. Let $K$ denote a maximal compact subgroup and
$T$ denote a maximal torus in $K$.  Let $h_\Z^\vee$ be the weight
lattice and $h_\Z$ be the co-weight lattice of $T$. Since $G$ is
simply connected, $Q^\vee \hookrightarrow h_\Z$ is an isomorphism.


Let $\F = \C((t))$ and $\O = \C[[t]]$.  The affine Grassmannian $\Gr
= \Gr_G$ is the ind-scheme $G(\F)/G(\O)$ (see~\cite{Kum} for more on
the ind-scheme structure).  It is a homogeneous space for the affine
Kac-Moody Group $\G$ associated to $W_\aff$, which for our purposes
can be taken simply to be $\G=G(\F)$.  The space $\Gr$ is
homotopy-equivalent to the space $\Omega K$ of based loops in $K$;
see~\cite{PS}.  Since $K$ is simply-connected, $\Gr$ is connected.
Each co-character $\lambda: S^1 \to T \,\in Q^\vee$ gives a point in
$\Omega K$, which we denote $t_\lambda$.  These are the $T$-fixed
points of $\Omega K$. We let $S = S(h^*_\Z) = H^T({\rm pt})$ denote
the symmetric algebra of $h^*_\Z$.  Since $T$ does not the contain
the rotation action of $S^1$ on $\Omega K$, the image of the
imaginary root $\delta = \alpha_0 + \theta$ in $h_\Z^*$ is 0.  In
other words, the map $Q_\aff \to h_\Z^*$ is given by $\alpha_0
\mapsto -\theta$ where $\theta$ is the highest root of the root
system $R$.

The group $\G$ possesses a {\it Bruhat decomposition} $\G =
\bigcup_{w \in W_\aff}\BB w \BB$ where $\BB$ denotes the Iwahori
subgroup. The Bruhat decomposition induces a decomposition of $\Gr$
into {\it Schubert cells} $\Omega_w = \BB w G(\O) \subset
G(\F)/G(\O)$:
\[
\Gr = \bigsqcup_{w \in W^0} \Omega_w = \bigcup_{w \in W^0} X_w,
\]
where the {\it Schubert varieties} $X_w$ are the closures of
$\Omega_w$.  In this paper we will be concerned with the homology
$H_*(\Gr)$, cohomology $H^*(\Gr)$, torus-equivariant homology
$H_T(\Gr)$ and torus-equivariant cohomology $H^T(\Gr)$ of the affine
Grassmannian. Note that the torus $T$ acts on $\Omega K$ by
(point-wise) conjugation. We will denote the Schubert classes in
homology, cohomology, equivariant homology and equivariant
cohomology as follows (see~\cite{Kum} for relevant definitions)
\[
\sigma_w \in H_*(\Gr) \;,\; \sigma^w \in H^*(\Gr) \;,\; \sigma_{(w)}
\in H_T(\Gr) \;,\; \sigma^{(w)} \in H^T(\Gr).
\]
More generally, we will use the notations $$\sigma^\p_w \in
H_*(\G/\p), \sigma^\p_{(w)} \in H_T(\G/\p), \sigma_\p^w \in
H^*(\G/\p), \sigma_\p^{(w)} \in H^T(\G/\p)$$ for the Schubert
classes for $\G/\p$ where $\p \subset \G$ is a parabolic subgroup.
Throughout this paper, all homology and cohomology rings will be
with $\Z$-coefficients.  We should note that the equivariant
homology theory that we use is equivariant Borel-Moore homology.
Since the stratification by Schubert cells is equivariant formal for
the $T$-action, we may think of  $H_T(\Gr)$ as the subspace of
$\Hom_S(H^T(\Gr),S)$ with basis $\sigma_{(w)}$ dual to
$\sigma^{(w)}$; see for example~\cite{Gra} (in particular Section
4.3) for the construction of these classes, which we omit. Both
$H_T(\Gr)$ and $H^T(\Gr)$ are free modules over $S = H^T({\rm pt})$
with basis the Schubert classes.

\smallskip
 The main aim of this
article is the identification of the classes $\sigma_w$ and
$\sigma^w$ explicitly, as polynomials.  These polynomials can be
considered affine analogues of Lascoux-Sc\"{u}tzenberger's Schubert
polynomials~\cite{LS}, which are polynomial representatives for
cohomology Schubert classes of the finite flag variety $K/T$.  Note
that there does not seem to be a classical, finite analogue of the
homology representatives $\sigma_w$, though in special cases the
homology polynomial representatives turn out to be Schur functions.



\section{Kostant-Kumar's NilHecke Ring}
\label{sec:nilhecke} The discussion of this section works for
symmetrizable Kac-Moody groups, but we will restrict our notation to
the affine case and make a small modification to the general theory
(Remark~\ref{rem:Kum}).  See~\cite{Kum, Pet} for further details.

\subsection{The affine nilHecke ring}
Let $\A_\aff$ denote the {\it affine nilHecke ring} of Kostant and
Kumar. (Note that Kostant and Kumar define the nilHecke ring over
the rationals, but we have found it more convenient, following
Peterson~\cite{Pet}, to work over $\Z$.)  It is the ring with a $1$
given by generators $\{A_i \mid i \in I \cup \{0\} \} \cup \{\lambda
\mid \lambda \in h_\Z^*\}$ and the relations
\begin{align*}
A_i \,\lambda &= (r_i \cdot \lambda)\, A_i +
\ip{\lambda,\alpha_i^\vee}\cdot1 & \mbox{for $\lambda \in h_\Z^*$,} \\
A_i\, A_i &= 0, \\
(A_iA_j)^m &= (A_jA_i)^m & \mbox{if $(r_ir_j)^m = (r_jr_i)^m$.}
\end{align*}
where the ``scalars'' $\lambda \in h_\Z^*$ commute with other
scalars.  The finite nilHecke ring is the subring of $\A_\aff$
generated by $\{A_i \mid i \in I \} \cup \{\lambda \mid \lambda \in
h_\Z^*\}$.  Combinatorially, one should think of the elements $A_i$
as {\it divided difference operators}.


Let $w \in W_\aff$ and let $w = r_{i_1} \cdots r_{i_l}$ be a reduced
decomposition of $w$.  Then $A_w := A_{i_1} \cdots A_{i_l}$ is a
well defined element of $\A_\aff$, where $A_\id = 1$.  By~\cite{KK}
or~\cite[Proposition 2-7]{Pet}, $\{A_w \mid w \in W_\aff\}$ is an
$S$-basis of $\A_\aff$.  The map $W_\aff \mapsto \A_\aff$ given by
$r_i \mapsto 1 - \alpha_i A_i \in \A_\aff$ is a homomorphism.
 Abusing notation, we write $w \in \A_\aff$ for the element in the
nilHecke ring corresponding to $w \in W_\aff$ under this map.  Then
$W_\aff$ is a basis for $\A_\aff$ over $\Frac(S)$ (not over $S$
since $A_i = \frac{1}{\alpha_i}(1-r_i)$).

Let $\A_0 \subset \A_\aff$ denote the subring over $\Z$ of $\A_\aff$
generated by the $A_i$ only, which we call the {\it affine
nilCoxeter algebra}.  There is a specialization map $\phi_0: \A_\aff
\rightarrow \A_0$ given by
$$
\phi_0: \sum_w a_w A_w \longmapsto \sum_w \phi_0(a_w) A_w
$$
where $\phi_0$ evaluates a polynomial $s \in S$ at 0.

For later use, we note the following straightforward result, whose
proof we omit; see~\cite[Proposition 4.30]{KK}.
\begin{lemma}
\label{lem:chev} Let $w \in W_\aff$ and $\lambda \in S$ be of degree
1. Then
\[
A_w \lambda = (w \cdot \lambda)A_w + \sum_{r_\alpha w \lessdot
w}\ip{\lambda, \alpha^\vee} A_{r_\alpha w}.
\]
Here $\lessdot$ denotes a cover in strong Bruhat order.
\end{lemma}
The coefficients $\ip{\lambda, \alpha^\vee}$ are known as {\it
Chevalley coefficients}.

\begin{remark}
\label{rem:Kum} The affine nilHecke ring $\A_\aff$ defined above is
slightly different to the nilHecke ring $\A(\hat{R})$ associated to
the affine root system $\hat{R}$ of $W_\aff$, as defined
in~\cite{KK}. Most importantly in our setup, the image of the
imaginary root $\delta = \alpha_0 + \theta$ in $h^*_\Z$ is 0.  In
other words, we do not consider the additional $S^1$-equivariance
obtained by rotating the loops in $\Omega K$. With this additional
equivariance, the equivariant homology $H_{T \times S^1}(\Gr)$ is
actually non-commutative (compare with~\cite{BFM}).  Another
difference is that the center of $\A(\hat{R})$ is much smaller than
that of $\A_\aff$. Shrawan Kumar has pointed out to us that the
center of $\A(\hat{R})$ is just the Weyl-group invariants of the
scalars.  We shall later find explicitly non-trivial elements in the
center of $\A_\aff$; see Section~\ref{sec:final}.
\end{remark}

\subsection{Action on the cohomology ring $H^T(\G/\p)$}
The ring $\A_\aff$ acts as generalized BGG-Demazure operators on
$H^T(X)$ for any $LK$-space $X$ (here $LK$ is the space of all
continuous loops into $K$).  The element $A_i$ corresponds to the
map $H^*(\G/\BB) \to H^{*-2}(\G/\BB)$ obtained by integration along
the fibers of the ${\mathbb P}^1$-fibration $\G/\BB \to \G/\p_i$
where $\p_i$ are the minimal parabolic subgroups; see~\cite{BGG,
Kum}.

Fix $\p$ a parabolic subgroup of $\G$ and let $W_\aff^P$ be the
minimal length coset representatives for the quotient of $W_\aff$ by
the corresponding parabolic subgroup $W_P \subset W_\aff$. Then in
terms of the Schubert classes $\sigma_\p^{(w)} \in H^T(\G/\p)$ and
$\sigma^\p_{(w)} \in H_T(\G/\p)$, where $w \in W_\aff^P$ we have
\begin{equation}
\label{eq:nilHeckeoncohom} A_v \cdot \sigma_\p^{(w)} = \begin{cases}
\sign(v)\sigma_\p^{(vw)} &
\mbox{if $\ell(v^{-1}) + \ell(vw) = \ell(w)$ and $vw \in W_\aff^\p$,}  \\
0 & \mbox{otherwise;} \end{cases}
\end{equation}
and
\begin{equation}
\label{eq:nilHeckeonhom} A_v \cdot \sigma_{(w)}^\p = \begin{cases}
\sigma^\p_{(vw)} &
\mbox{if $\ell(vw) = \ell(v)+ \ell(w)$ and $vw \in W_\aff^P$,}  \\
0 & \mbox{otherwise.} \end{cases}
\end{equation}
Note that the formulae above may differ from some of the literature
by either a sign or the replacement of left multiplication by right
multiplication.
\medskip

Let $w \in W_\aff$ and consider $\iota_w : {\rm pt} \to \G/\p$ given
by $* \mapsto w\p$, the inclusion of a $T$-fixed point. Then the
localization maps $\psi_w = \iota_w^*: H^T(\G/p) \rightarrow
H^T({\rm pt}) = S$ form a basis for $H_T(\G/\p)$ over $\Frac(S)$.
Note that we have $\psi_w = \psi_{w'}$ if $w \in w' W_P$.  The
change of basis matrix between $\{\psi_w\}$ and
$\{\sigma_{(w)}^\p\}$ is given by the $D$-matrix of~\cite{KK}.  The
action of the nilHecke ring on the basis $\{\psi_w\}$  is given by
\begin{equation}
\label{eq:nilHeckeFixed}
 v \cdot \psi_w = \psi_{vw}.
\end{equation}

\begin{remark}
In~\cite{Pet}, the nilHecke ring $\A_\aff$ is identified with
$H_T(\G/\BB)$ via the map
$$
a \in \A_\aff \longmapsto \psi_\id \circ a \in
\Hom_S(H^T(\G/\BB),S).
$$
One can give $H^T(\G/\BB)$ a product and coproduct which agrees with
that of $\A_\aff$ (see Section~\ref{sec:coproduct}).
\end{remark}

\subsection{The coproduct on $\A_\aff$}
\label{sec:coproduct} Define the coproduct map $\Delta: \A_\aff
\rightarrow \A_\aff \otimes_S \A_\aff$ by
\begin{align*}
\Delta(s) &= 1 \otimes s = s \otimes 1 & \mbox{for $s \in S$} \\
\Delta(A_i) &= A_i \otimes 1 + r_i \otimes A_i = 1 \otimes A_i + A_i
\otimes r_i \\ &= A_i \otimes 1 + 1 \otimes A_i - A_i \otimes
\alpha_iA_i.
\end{align*}
This is a well defined map, which in addition is cocommutative. One
can deduce from these relations that $\Delta(w) = w \otimes w$. (In
the original work of~\cite{KK}, this last relation was used to
define $\Delta$, but we shall follow the set up of~\cite{Pet}).

One should be careful since the tensor product $\A_\aff \otimes_S
\A_\aff$ is {\it not} a ring in the obvious way. For example,
\[
(A_i \otimes 1 ).(1 \otimes \alpha_i) \neq (A_i \otimes 1
).(\alpha_i \otimes 1)
\]
However, it is shown in~\cite{Pet} that the above formulae still
give an action of $\A_\aff$ on $\A_\aff \otimes_S \A_\aff$. That is,
$\Delta(a) = a \cdot (1 \otimes 1)$ for any $a \in \A$.  In
particular, if $M$ and $N$ are $\A_\aff$ modules, then so is $M
\otimes_S N$.

Note that $\phi_0$ also sends $\A_\aff \otimes_S \A_\aff$ to $\A_0
\otimes_\Z \A_0$ by evaluating the coefficients at 0 when writing in
the basis $\{A_w \otimes A_v\}_{w,v \in W_\aff}$.  The following
theorem relates the coproduct on $\A_\aff$ with the equivariant
cohomology $H^T(\G/\BB)$.

\begin{theorem}[\cite{KK,Ara}]
\label{thm:kk} Let $$ \Delta(A_w) = \sum_{u,v \in W_\aff} a_w^{u,v}
A_u \otimes A_v.$$  Then $a^{u,v}_w$ are the (Schubert) structure
constants of $H^T(\G/\BB)$, so that $$\sigma_\BB^{(u)}\cdot
\sigma_\BB^{(v)} = \sum_{w \in W_\aff} a_w^{u,v} \sigma_B^{(w)},$$
where $\sigma_\BB^{(w)}$ denote the cohomology Schubert classes of
$\G/\BB$.
\end{theorem}
Note that since $H^T(\G/\p)$ imbeds into $H^T(\G/\BB)$, the
structure constants for $\G/\BB$ include the structure constants for
$\G/\p$.  Theorem~\ref{thm:kk} is in fact valid for all
symmetrizable Kac-Moody groups.


\section{Equivariant homology of $\Gr$ and Peterson's $j$-homomorphism}
\label{sec:pet} Peterson~\cite{Pet} has constructed an isomorphism,
called $j$, between the equivariant homology $H_T(\Gr)$ of the
affine Grassmannian and the centralizer $Z_{\A_\aff}(S)$ of $S$ in
$\A_\aff$.  We shall call the subalgebra $Z_{\A_\aff}(S)$ the {\it
Peterson subalgebra}.  Since Peterson's work~\cite{Pet} has yet to
be published, we include here a simplified construction of this
isomorphism.  Our definition of the isomorphism is more axiomatic
than Peterson's, at the expense of being less natural; see
Remark~\ref{rem:pet}.

\smallskip

For this section we will think of $\Gr$ both as the based loops
$\Omega K$ and as $\G/\p$ for a maximal parabolic $\p$, but will use
the former notation.  Note that $W_P = W$ in this case.  Recall that
$H_T(\Omega K)$ is the $S$-submodule of $\Hom_S(H^T(\Omega K),S)$
spanned by the classes $\sigma_{(w)}$. Over ${\rm Frac}(S)$,
$H_T(\Omega K)$ is also spanned by the classes $\{\psi_t \mid t =
t_\lambda \in Q^\vee \subset W_\aff \}$.


Since $\Omega K$ is a group and the multiplication $\Omega K \times
\Omega K \rightarrow \Omega K$ is $T$-equivariant, $H_T(\Omega K)$
and $H^T(\Omega K)$ obtain the structures of dual Hopf-algebras. The
Hopf-structure maps of $H_T(\Omega K)$ can be calculated directly in
terms of the $\psi_t$ as follows ($\counit$ denotes the counit and
$c$ denotes the antipode):
\begin{align}
\label{eq:omega} \psi_{\rm id} &= 1, & \counit(\psi_t) &=1,  &
c(\psi_t) &=\psi_{t^{-1}}, \\
\Delta \psi_t &= \psi_t \otimes \psi_t, &\psi_t \psi_{t'} &=
\psi_{tt'}.
\end{align}
The scalars $s \in S \subset H_T(\Omega K)$ are central and $c(s) =
s$.  In particular, $H_T(\Omega K)$ is a commutative and
co-commutative Hopf-algebra; see also~\cite{BFM} for a comparison.

\begin{lemma}
\label{lem:commute} Let $t_\lambda \in \A_\aff$ be the image of a
translation $t_\lambda \in W_\aff$ in the nilHecke ring.  Then
$t_\lambda \in Z_{\A_\aff}(S)$.
\end{lemma}
\begin{proof}
The translations $t_\lambda$ act on $Q_\aff = \oplus_{i \in I \cup
\{0\}} \alpha_i$ via $t_\lambda \cdot (\alpha + n \delta) = \alpha +
n \delta - \ip{\alpha, \lambda}\delta$, where $\alpha\in R$ is a
root of the finite root system. Since the image of $\delta$ in $S$
is 0, we see that $t_\lambda$ lies in $Z_{\A_\aff}(S)$.
\end{proof}

\begin{definition}
\label{def:j} Define a map $j: H_T(\Omega K) \to \A_\aff$ by
$$
\psi_t \longmapsto t
$$
for all $t \in Q^\vee$ and extending by linearity (over $S$).
\end{definition}

From the definition, the image of $j$ lies in  $\Frac(S)\otimes_S
\A_\aff$ but it turns out that the image of $j$ lies in $(\Frac(S)
\otimes {\rm span}_S\{t \mid t \in Q^\vee\}) \cap \A_\aff$.  This
will be clear after Lemma~\ref{lem:action}.


Thus $j(H_T(\Omega K))$ is the subalgebra of $\A_\aff$ generated by
the elements $t \in Q^\vee$ and $s \in S$.  Note that by
Lemma~\ref{lem:commute}, $t \in Z_{\A_\aff}(S)$ so that the image of
$j$ lies in $Z_{\A_\aff}(S)$ and in particular $j(H_T(\Omega K))$ is
an algebra over $S$.

\begin{lemma}
\label{lem:action} The map $j$ is a homomorphism of rings.  It is
compatible with the $\A_\aff$ action on $H_T(\Omega K)$ as follows:
$$
\sigma \, \sigma' = j(\sigma) \cdot \sigma',
$$
for $\sigma, \sigma' \in H_T(\Omega K)$.
\end{lemma}
\begin{proof}
Since both $H_T(\Omega K)$ and $j(H_T(\Omega K))$ are algebras over
$S$ (that is, $S$ is central), we need only compute on the basis
$\psi_t$. The first statement follows from $\psi_t \psi_{t'} =
\psi_{tt'}$. The second statement follows from formula
(\ref{eq:nilHeckeFixed}).
\end{proof}

Every $A_w$ for $w \in W_\aff$ acts non-trivially on $H_T(\Omega
K)$.  By writing $j(\sigma) = \sum_{w \in A_w} a_\sigma^w \,A_w$ and
using (\ref{eq:nilHeckeoncohom}), Lemma~\ref{lem:action} shows that
the coefficients $a_\sigma^w$ lie in $S$, not $\Frac(S)$.


\begin{theorem}[\cite{Pet}]
\label{thm:pet} There is an isomorphism $j: H_T(\Omega K)
\rightarrow Z_{\A_\aff}(S)$ such that
\[
j(\sigma_{(x)}) = A_x  \mod I
\]
where $x \in W^0$ and
\[
I = \sum_{w \in W \,;\; w \neq \id} \A \cdot A_w.
\]
\end{theorem}
\begin{proof}
We have already shown that $j(H_T(\Omega K)) \subset
Z_{\A_\aff}(S)$.  By Lemma~\ref{lem:action}, we have
$j(\sigma_{(x)}).\sigma_{(\id)} = \sigma_{(x)}$ for $x \in W^0$.
Using (\ref{eq:nilHeckeonhom}), we see that $j(\sigma_{(x)}) = A_x +
a$ where $a$ lies in the annihilator of $\sigma_{(\id)} \in
H_T(\Omega K)$.  This annihilator is easily seen to be equal to the
left ideal $I$ defined above.


To show that $j$ is an isomorphism, write $a \in Z_{\A_\aff}(S)$ as
$$a = \sum_{w \in W_\aff}a_w \, w.$$  Now we use the commutation
formula $w \, s = (w \cdot s)\, w$ in $\A_\aff$ (see~\cite{KK}).
Since $W$ acts faithfully on $S$, if $w \notin Q^\vee$ we must have
$a_w =0$, and so $a \in j(H_T(\Omega K))$.
\end{proof}

%
%
%

The isomorphism of Theorem~\ref{thm:pet} can be made into a
Hopf-isomorphism by giving $Z_{\A_\aff}(S)$ the restriction of the
coproduct $\Delta$ on $\A$ to $Z_{\A_\aff}(S)$.  Since we have
$\Delta(t) = t \otimes t \in Z_{\A_\aff}(S) \otimes_S
Z_{\A_\aff}(S)$ this agrees with the coproduct of $H_T(\Omega K)$.
We define the antipode by $c(t) = t^{-1}$ for $t \in Q^\vee \subset
Z_{\A_\aff}(S)$ (see (\ref{eq:omega})). It is induced by $\alpha_i
\mapsto -\alpha_{\omega(i)}$ for $i \in I$ ($\alpha_0 \mapsto
-\alpha_0$) and $A_i \mapsto -A_{\omega(i)}$ where $\omega$ is the
diagram automorphism obtained by conjugation by the longest element
of $W$.

The elements $j(\sigma_{(t_\lambda)})$ for antidominant coweights
$\lambda$ can be calculated explicitly.

\begin{prop}
\label{prop:specialclass} Let $x \in W^0$ and $\lambda \in Q^-$.
Then we have
\begin{equation}
\label{eq:hommult} \sigma_{(x)}\sigma_{(t_\lambda)} =
\sigma_{(xt_\lambda)}
\end{equation}
in $H_T(\Omega K)$.  Explicitly one has
\begin{equation}
\label{eq:conjcent} j(\sigma_{(t_\lambda)}) = \sum_{w \in
W/W_\lambda} A_{t_{w \cdot \lambda}}
\end{equation}
where $W_\lambda \subset W$ is the stabiliser of $\lambda$.
\end{prop}
\begin{proof}
Write $t := t_\lambda$.  First we note that $A_i \cdot \sigma_{(t)}
= 0$ for $i \neq 0$. This follows from (\ref{eq:nilHeckeonhom}),
since using (\ref{eq:length}) one can check that for each $i \neq 0$
we have either $\ell(r_i t) < \ell(t)$ or $r_i t \notin W^0$.  Now
we compute, using Theorem~\ref{thm:pet}
\begin{align*}
\sigma_{(x)}\sigma_{(t)} &=  j(\sigma_{(x)}) \cdot \sigma_{(t)} \\
&= (A_x + a) \cdot \sigma_{(t)} & \mbox{where $a \in I$.}\\
&= A_x \cdot \sigma_{(t)} & \mbox{since $a \cdot \sigma_{(t)} = 0$.} \\
&= \sigma_{(xt)}.
\end{align*}
For the last statement, we write, using the fact that $H_T(\Omega
K)$ is commutative, $\sigma_{(xt)} = j(\sigma_{(t)}) \cdot
\sigma_{(x)}$.  Now let $x = w t'$ where $w \in W$ and $t' \in
Q^\vee$.  Then $$xt = wt't_\lambda = wt_\lambda t' = t_{w\cdot
\lambda}w t'$$ so that $A_{t_{w\cdot \lambda}}$ must occur in the
expansion of $j(\sigma_{(t)})$.  But the element $ \sum_{w \in
W/W_\lambda} A_{t_{w \cdot \lambda}}$ lies in $Z_{\A_\aff}(S)$ and
is of the form $(A_t \mod I)$, so by Theorem~\ref{thm:pet} one
obtains (\ref{eq:conjcent}).
\end{proof}

\smallskip

\begin{remark}
\label{rem:pet} 

In~\cite{Pet}, the map $j$ is obtained in the context of ``compact
characteristic operators'' by considering the action
$$
\phi: \Omega K \times LK/T \longrightarrow LK/T
$$
of $\Omega K$ on $\G/\BB \simeq LK/T$ where $LK$ is the space of all
continuous loops into $K$.  This action is obtained via the
inclusion $\Omega K \hookrightarrow LK$ and is $T$-equivariant.
Each $\sigma \in H_T(\Omega K)$ induces a composition
$$
 H^T(LK/T) \longrightarrow H^T(\Omega K) \otimes_S
H^T(LK/T) \longrightarrow S \otimes_S H^T(LK/T)$$ where the first
map is $\phi^*$ and the second map is $\sigma \otimes \id$.  This
defines a $H_T(\Omega K)$-module structure on $H^T(LK/T)$.  If
$\sigma \in H_T(\Omega K)$ then we obtain an operator $a_\sigma$ on
$H^T(LK/T)$.  This operator $a_\sigma$ acts on $H^T(LK/T)$ in the
same way as $j(\sigma) \in \A_\aff$.  The expansion of $j(\sigma)$
in the $\{A_w \mid w \in W_\aff\}$ basis can be obtained by
calculating the action of $a_\sigma$ on the Schubert basis $\{
\sigma^{(w)}\}$, using (\ref{eq:nilHeckeoncohom}).

Note that since $\Omega K \hookrightarrow LK/T$ takes a fixed point
$t_\lambda \in \Omega K$ to the fixed point $t_\lambda \, T$, it is
immediate that this definition agrees with Definition~\ref{def:j}.
\end{remark}

%
%

\section{Homology of affine Grassmannian and Fomin-Stanley subalgebra}
\label{sec:FShom}
\subsection{An identity for finite Weyl groups} In this subsection we
let $W$ be a finite Weyl group and $H^*(K/T)$ be the cohomology of
the corresponding flag variety.  Also let $w^\circ$ denote the
longest element of $W$.

\begin{proposition}
\label{prop:finite} Suppose that for some coefficients $\{b_u \in
\Z\}_{u \in W}$ the following identity holds in $\Z[W]$ for all
integral weights $\lambda \in h^*_\Z$:
\[
\sum_{u \in W; \,\,l(u) > 0} b_{u} \sum_{ u r_\alpha  \lessdot u}
\ip{\lambda,\alpha^\vee} u r_\alpha = 0.
\]
Then $b_{u} = 0$ for all $u$.
\end{proposition}
\begin{proof}
First apply the transformation $u \mapsto w^\circ u$ to the identity
of the Proposition.  Then reindexing the $b_u$, we obtain
\[
\sum_{u \in W; \,\,u \neq w^\circ} b_{u} \sum_{ u r_\alpha \gtrdot
u} \ip{\lambda,\alpha^\vee} u r_\alpha = 0
\]
for all $\lambda$.

Let $\zeta^u \in H^*(K/T)$ denote the cohomology Schubert classes of
the finite flag variety. By the Chevalley-Monk formula~\cite{BGG} we
have
\[
[\lambda] \cdot \zeta^u = \sum_{u r_\alpha \gtrdot u}
\ip{\lambda,\alpha^\vee}\zeta^{u r_\alpha}
\]
where $[\lambda] \in H^*(K/T)$ denotes the image of $\lambda  \in
h^*_\Z$ under the characteristic homomorphism $S(h^*_\Z) \rightarrow
H^*(K/T)$. For example, if $\lambda = \omega_i$ is a fundamental
weight then $[\omega_i] = \zeta^{s_i}$.  It is well known that $
\zeta^{s_1}, \zeta^{s_2},\ldots, \zeta^{s_{n-1}}$ generate the
positive degree part of $H^*(K/T)$ or alternatively that the
characteristic homomorphism is surjective.

Suppose that $[\lambda] \cdot \zeta= 0$ for some $\zeta \in
H^*(K/T)$ and all $\lambda \in h^*_\Z$.   Then $\zeta^u\cdot \zeta =
0$ for all $u \neq \id$.  If $\ell(v) + \ell(u) = \ell(w_\circ)$ we
have $\zeta^{v} \cdot \zeta^u = \delta_{v,w_\circ
u}\zeta^{w_\circ}$.  Thus we find that $\zeta$ must be a multiple of
the class $\zeta^{w_\circ}$.  Letting $\sigma= \sum_u b_u \zeta^u$
and applying the Chevalley-Monk formula we obtain the proposition.
\end{proof}

\subsection{Affine Fomin-Stanley subalgebra}
\label{sec:FS} Define the {\it affine Fomin Stanley subalgebra} as
the subspace $\B' \subset \A_0$ satisfying:
\[
\B' = \{ a \in \A_0 \mid \phi_0(as) = \phi_0(s)a \,\, \text{for all
$s \in S$} \}.
\]
The name is justified by the following computation.
\begin{lemma}
The space $\B'$ defined above is a subalgebra.
\end{lemma}
\begin{proof}
Let $a,b \in \B'$ and $s \in S$.  Then $bs = \phi_0(s)b + \sum_{w
\in W_\aff}c_w A_w$ where $c_w \in S$ and $\phi_0(c_w) = 0$.  Thus
$$
(ab)s = a \phi_0(s) b+ \sum_{w \in W_\aff} \phi_0(c_w) a A_w =
\phi_0(s)(ab),
$$
and hence $ab \in \B'$.
\end{proof}

A combinatorial definition of the same algebra, in type $A$, will be
given in Section~\ref{sec:kschur}.


\begin{proposition}
\label{prop:grass} Let $b \neq 0 \in \B'$ and write $b = \sum_{w \in
W_\aff} b_w A_w$ with $b_w \in \Z$. Then $b_w \neq 0$ for some $w
\in W^0$.
\end{proposition}

\begin{proof}
Let $D = \{w \in W_\aff \mid b_w \neq 0\}$.  For each $w \in W_\aff$
we may uniquely write $w = x_wy_w$ where $x_w \in W^0$ and $y_w \in
W$. Let $d = \{\min(\ell(y_w)) \mid w \in D\}$.  We write
$\ell_0(w):=\ell(y_w)$.

Suppose $d \neq 0$ and let $w \in D$ minimize $\ell_0(w)$.  Let
$\lambda \in S$ be of degree 1.  Then by Lemma~\ref{lem:chev},
$\phi_0(A_w \lambda) =  \sum_{ w r_\alpha \lessdot w}\ip{\lambda,
\alpha^\vee} A_{w r_\alpha}$.  We know that $w \gtrdot v$ if and
only if a reduced decomposition of $v$ is obtained from a reduced
decomposition of $w$ by removing a simple generator. Since $w =
x_wy_w$, each such $v$ satisfies $\ell_0(v) \geq \ell_0(w) - 1$. Let
$D_w = \{v \lessdot w \mid \ell_0(v) = \ell_0(w) - 1\}$.  Then $v
\in D_w$ if and only if $v = x_vy_v$ where $x_v = x_w$ and $y_v
\lessdot y_w$.  In this case $v = wr_\alpha$, where $\alpha \in R$
is a root of the finite root system.

Now write $\phi_0(b \lambda) = \sum_v b'_v A_v$ and focus only on
the coefficients of $b'_v$ satisfying $\ell_0(v) = d - 1$ and $v =
xy_v$ for some fixed $x \in W^0$.  If $b \in \B'$ then $b'_v = 0$.
Thus in particular, for every $\lambda \in S$ of degree 1, we have
\[
\sum_{u \in W \,\mid \,\ell(u) = d} b_{x u} \sum_{ u r_\alpha
\lessdot u} \ip{\lambda,\alpha^\vee} A_{x}A_{u r_\alpha} = 0.
\]
Now the set $\{A_xA_u \mid u \in W\}$ is independent over $S$.  So
by replacing $A_xA_u$ by $u \in \Z[W]$, we see that this is
impossible by Proposition~\ref{prop:finite}.  Thus we conclude that
we must have $d = 0$.
\end{proof}

The following theorem says that evaluation at 0 maps the Peterson
subalgebra onto the affine Fomin Stanley subalgebra.

\begin{prop}
\label{prop:bprime} We have $\phi_0(Z_{\A_\aff}(S)) = \B'$.  More
precisely, the set $\{\phi_0(j(\sigma_{(u)}))\mid u \in W^0\}$ forms
a basis of $\B'$ over $\Z$.  Furthermore, the element
$\phi_0(j(\sigma_{(u)}))$ is the unique element in $\B'$ with unique
Grassmannian term $A_u$.
\end{prop}

\begin{proof} The fact that $\phi_0(Z_{\A_\aff}(S)) \subset \B'$ is a trivial
calculation. Now let $b \in \B'$.  By Proposition~\ref{prop:grass}
it contains a Grassmannian term $A_u$ with non-zero coefficient
$b_u$.  By Theorem~\ref{thm:pet}, $b - b_u\phi_0(j(\sigma_{(u)}))$
has strictly fewer Grassmannian terms and also lies in $\B'$.
Repeating, we see that one can write $b$ uniquely as a $\Z$-linear
combination of the elements $\phi_0(j(\sigma_{(u)}))$.  The last
statement follows from Theorem~\ref{thm:pet}.
\end{proof}

By Theorem~\ref{thm:pet}, we have $\Delta(z) \in Z_{\A_\aff}(S)$ for
$z \in Z_{\A_\aff}(S)$ where $\Delta$ is the coproduct of $\A_\aff$
introduced in Section~\ref{sec:coproduct}.  Applying the evaluation
$\phi_0$ (and using Theorem~\ref{thm:kk}) we see that $\B'$ is
closed under the operation $\Delta_{\B'} = \phi_0 \circ \Delta$.
Thus $\B'$ attains the structure of a Hopf algebra over $\Z$.


Since we have $H^*(\Gr) = H^T(\Gr) \otimes_S \Z$, where $S$ acts on
$\Z$ by evaluation at 0, the product and coproduct structure
constants of $H^*(\Gr)$ (respectively $H_*(\Gr)$) are obtained from
those of $H^T(\Gr)$ (respectively $H_T(\Gr)$) by evaluating at 0.
We thus obtain the following theorem.

\begin{theorem}
\label{thm:bprime} The affine Fomin-Stanley subalgebra $\B'$ is a
model for the homology $H_*(\Gr)$.  More precisely, the map
$H_*(\Gr) \to \B'$ given by $\sigma_{u} \mapsto
\phi_0(j(\sigma_{(u)}))$ is an isomorphism of Hopf algebras.
\end{theorem}

As a consequence of the commutativity of $H_*(\Gr)$ we have the
following result.
\begin{corollary}
The affine Fomin-Stanley subalgebra $\B'$ is a commutative algebra.
\end{corollary}

In the forthcoming sections, our aim will be to explicitly describe
the elements $\phi_0(j(\sigma_{(u)}))$ in the type $A_{n-1}$ case;
that is, with $K = SU(n)$.

\begin{remark}
As will be clear shortly, it is reasonable to call the homology
Schubert basis elements $\phi_0(j(\sigma_{(u)}))$ ``non-commutative
$k$-Schur functions'' for the group $G$.  The element
$\phi_0(j(\sigma_{(u)}))$ is characterized by being an element of
$\B'$ and having a unique Grassmannian term $A_u$.
\end{remark}


\section{$k$-Schur functions}
\label{sec:symkschur} For the remainder of the paper, unless
specified otherwise, we will restrict ourselves to type $A_{n-1}$.
We have $G = SL(n,\C)$, $K = SU(n)$, $\G = SL(n,\F)$.  Thus $W =
S_n$ is the symmetric group, $W_\aff$ is the affine symmetric group
and the affine simple roots are given by $\{\alpha_i \mid i \in
\Z/n\Z\}$.

\subsection{Symmetric functions} \label{sec:sym} We refer
to~\cite{Mac} for details concerning the material of this section.
Let $\Lambda = \Lambda_\Z$ denote the ring of symmetric functions
over $\Z$ in infinitely many variables $x_1,x_2,\ldots$. We write
$h_i(x)$ for the {\it homogeneous symmetric functions} and for a
partition $\lambda = (\lambda_1 \geq \lambda_2 \geq \cdots)$, we
write $h_\lambda(x) = h_{\lambda_1(x)}h_{\lambda_2(x)}\cdots$. The
elements $h_1(x),h_2(x),\ldots \in \Lambda$ form a set of
algebraically independent set of generators of $\Lambda$. Similarly,
we let $e_i(x)$ denote the {\it elementary symmetric functions}. The
algebra involution $\omega$ of $\Lambda$ is defined by setting
$\omega(h_i(x)) = e_i(x)$.

We let $\{m_\lambda(x) \in \Lambda \mid \lambda \ \text{a
partition}\}$ denote the {\it monomial symmetric functions}. They
form a basis of the ring of symmetric functions over the integers.
Also let $p_k(x)$ denote the {\it power sum symmetric functions}.
The power sum symmetric functions are algebraically independent
generators of the ring of symmetric functions over $\Q$.

Let $\Lambda_n \subset \Lambda$ denote the subring of the symmetric
functions generated by $h_i(x)$ for $i \in [0,n-1]$.  Let
$\Lambda^n$ denote the quotient algebra of $\Lambda$ given by
$\Lambda^n = \Lambda/\langle m_\lambda(x) \mid \lambda_1 \geq n
\rangle$. Clearly the set $\{m_\lambda(x) \mid \lambda_1 < n\}$
forms a basis of $\Lambda^n$.  When giving an element $\bar f \in
\Lambda^n$ we will usually just give a representative $f \in
\Lambda$ without further comment.

The {\it Hall inner product}, denoted $\ip{.,.}: \Lambda \times
\Lambda \rightarrow \Z$, is a symmetric non-degenerate pairing
defined by $\ip{h_\lambda(x),m_\mu(x)} = \delta_{\lambda\mu}$. It
induces a non-degenerate pairing $\ip{.,.}: \Lambda_n \times
\Lambda^n \rightarrow \Z$.  The ring of symmetric functions
$\Lambda$ can be given the structure of a (commutative and
co-commutative) Hopf algebra with coproduct $\Delta_\sym:\Lambda(x)
\rightarrow \Lambda(x) \otimes \Lambda(y)$ given by
$\Delta_\sym(h_i(x))= \sum_{0 \leq j \leq i} h_j(x) \otimes h_i(y)$,
where $h_0 = 1$.  The coproduct can alternatively be given by
specifying $\Delta_\sym(p_i(x)) = p_i(x) \otimes 1 + 1 \otimes
p_i(y)$.  This Hopf-algebra structure gives $\Lambda_n$ and
$\Lambda^n$ the structures of dual Hopf algebras.  The counit
$\counit$ is given by the constant term $\counit(f(x)) = f(0)$. The
antipode $c$ is given by $c(h_i(x)) = (-1)^ie_i(x)$; thus if $f(x)
\in \Lambda$ is homogeneous of degree $d$ then $c(f(x)) = (-1)^d
\omega(f(x))$.

\smallskip

The homology and cohomology rings (and their Hopf-algebra
structures) of $\Omega K = \Omega SU(n)$ (in fact for the based
loops of all compact Lie groups) were earlier computed by Bott.

\begin{theorem}[\cite{Bot}]
\label{thm:bott} We have an isomorphism in homology
\[
H_*(\G/\p) = \Z[\sigma_1,\sigma_2,\ldots,\sigma_{n-1}]
\]
where ${\rm dim}\, \sigma_i = 2i$ and the coproduct is given by
$\Delta(\sigma_i) = \sum_{0 \leq j \leq i} \sigma_j \times
\sigma_{i-j}$ where $\sigma_0 = 1$. In cohomology we have
\[
H^*(\G/\p) = {\mathcal S}H^*(\mathbb{CP}^{n-1})
\]
where ${\mathcal S}$ denotes an infinite symmetric power.  The
primitive subspace of $H^*(\G/\p)$ is the space spanned by the power
sums $p_k(u) = u_1^k + u_2^k + \cdots  $ where the $u_i$ are the
generators of the different copies of $H^*(\mathbb{CP}^{n-1})\simeq
\Z[u]/u^{n}$.
\end{theorem}

Thus the homology $H_*(\G/\p)$ can be identified with the subring of
symmetric functions $\Lambda_n$ via $\sigma_i \mapsto h_i$ while the
cohomology $H^*(\G/\p)$ can similarly be identified with the
quotient $\Lambda^n$.  The aim of this paper is thus to identify the
Schubert classes $\sigma_w \in H_*(\G/\p) \simeq \Lambda_n$ and
$\sigma^w \in H^*(\G/\p) \simeq \Lambda^n$ as explicit symmetric
functions.

Note that because of the nature (one dimensional for each graded
piece with degree from $0$ to $n-1$) of the primitive subspace of
$H^*(\G/\p)$ there is essentially no choice in these
Hopf-isomorphisms.  The only flexibility would be sending $u_1^k +
u_2^k + \cdots $ to the negative power sum $-p_k$ instead of $p_k$.
This is ruled out if we want the Schubert classes $\sigma^w$ to be
represented by symmetric functions with a positive monomial
expansion.

\subsection{Affine Schur functions and $k$-Schur functions}
\label{sec:kschur}
 In this section we define dual bases
$\{s_\lambda^{(k)}(x)\}$ of $\Lambda_n$ and $\{F_\lambda(x)\}$ of
$\Lambda^n$ called respectively the {\it $k$-Schur functions}, and
the {\it affine Schur functions} or {\it dual $k$-Schur functions}.
Historically, the $k$-Schur functions $\{s_\lambda^{(k)}(x)\}$ were
introduced first in~\cite{LLM}, and were further studied
in~\cite{LM,LM04}.  Their introduction was motivated by a series of
conjectures the most striking of which is their relationship with
the Macdonald polynomials, which we briefly explain.

The original $k$-Schur functions $s^{(k)}_\lambda(x;t)$ depended on
an extra parameter $t$, and the $k$-Schur functions we use are the
$t=1$ specializations: $s^{(k)}(x) = s^{(k)}_\lambda(x;1)$.  Let
$H_\mu(x;q,t)$ be given by the plethystic substitution
$H_\mu(x;q,t)= J_\mu(x/(1-q);q,t)$ where $J_\mu(x;q,t)$ is the
integral form of Macdonald polynomials~\cite{Mac}.  Define
$K^{(k)}_{\nu\mu}(q,t) \in \Z[q,t]$ and $\pi^{(k)}_{\ll\nu}(t) \in
\Z[q,t]$ by
\[
H_\mu(x;q,t) = \sum_\nu K^{(k)}_{\nu\mu}(q,t) s^{(k)}_\nu(x;t)
\;\;\; ; \;\;\; s^{(k)}_\nu(x;t) = \sum_\ll \pi^{(k)}_{\ll\nu}(t)
s_\ll(x).
\]
Here $s_\ll(x)$ denotes a Schur function.  Then it is conjectured
that $K^{(k)}_{\nu\mu}(q,t) \in \N[q,t]$ and $\pi^{(k)}_{\ll\nu}(t)
\in \N[t]$ which would refine the (proven) ``Macdonald positivity
conjecture'' that the Schur expansion of $H_\mu(x;q,t)$ has
coefficients in $\N[q,t]$; see~\cite{Hai}.
\smallskip

There are a number of conjecturally equivalent definitions of
$k$-Schur functions.  We will use a reformulation of the definition
from~\cite{LM04} which involves first defining the {\it affine
Stanley symmetric functions} from~\cite{Lam}.

\begin{definition}[\cite{Lam}]
Let $a = a_1a_2\cdots a_k$ be a word with letters from $\Z/n\Z$ so
that $a_i \neq a_j$ for $i \neq j$.  Let $A = \{a_1,a_2,\ldots,
a_k\} \subset \Z/n\Z$.  The word $a$ is \emph{cyclically decreasing}
if for every $i$ such that $i,i+1 \in A$, the letter $i+1$ precedes
$i$ in $a$.  An affine permutation $w \in W_\aff$ is
\emph{cyclically decreasing} if $w = s_{a_1} \cdots s_{a_k}$ for
some cyclically decreasing sequence $a_1a_2\cdots a_k$.
\end{definition}

Define the elements $h_i \in \A_0 \subset \A_\aff: i \in [0,n-1]$ by
the formula
\[
h_i = \sum_w A_w
\]
where the sum is over all cyclically decreasing permutations $w \in
W_\aff$ with length $l(w) = i$.  If $I \subset \Z/n\Z$ and $w$ is
the corresponding cyclically decreasing permutation then we will
write $A_I$ for $A_w$.  This gives a one-to-one correspondence
between cyclically decreasing permutations and subsets of $\Z/n\Z$.
Here $h_0 = A_{\id} = 1$.  As an example, with $n = 4$ and $i = 2$
we have $h_3 = A_3A_2 + A_3A_1 + A_0A_3 + A_2A_1 + A_2A_0 + A_1A_0$.

\medskip

Let $\B$ denote the subalgebra of $\A_0 \subset \A$ generated by the
$h_i$ for $i \in [0,n-1]$, which we call the {\it affine
Fomin-Stanley subalgebra} (we will show that $\B$ agrees with $\B'$
defined in Section~\ref{sec:FS}).

\begin{theorem}[\cite{Lam}]
\label{thm:sym} The algebra $\B$ is commutative.  It is isomorphic
to the subalgebra $\Lambda_n$ of the symmetric functions generated
by the homogeneous symmetric functions $h_i(x)$ for $i \in [0,n-1]$,
under the map $\psi: h_i(x) \mapsto h_i \in \B$.
\end{theorem}

Let $\ip{.,.}: \A_0 \times \A_0 \rightarrow \Z$ denote the symmetric
non-degenerate pairing defined by $\ip{A_w,A_v} = \delta_{wv}$.

\begin{definition}[\cite{Lam}]
\label{def:affineStanley} Let $w \in W_\aff$.  Define the
\emph{affine Stanley symmetric functions} $\tF_w(x) \in \Lambda$ by
\[\tF_{w}(x) = \sum_{a = (a_1,a_2,\ldots,a_t)}
\ip{h_{a_t}h_{a_{t-1}} \cdots h_{a_1} \cdot 1, A_w} x_1^{a_1}
x_2^{a_2} \cdots x_t^{a_t},\] where the sum is over compositions of
$\ell(w)$ satisfying $a_i \in [0,n-1]$.
\end{definition}

This definition is similar to the definition of Stanley symmetric
functions given by Fomin and Stanley~\cite{FS}.

The image in $\Lambda^n$ of the set $\{\tF_w(x) \mid w \in W^0\}$
forms a basis of $\Lambda^n$ (see~\cite{Lam, LM05}).  We called
these functions {\it affine Schur functions} in~\cite{Lam}.  They
were earlier introduced in a different manner in~\cite{LM05}, where
they were called {\it dual $k$-Schur functions}.

Define the {\it $k$-Schur functions} $\{s_w^{(k)}(x) \mid w \in
W^0\}$ as the dual basis of $\Lambda_n$ to the affine Schur
functions under the Hall inner product, where $k = n-1$.  There is a
bijection $w \leftrightarrow \lambda(w)$ from affine Grassmannian
permutations $\{w \in W^0\}$ to partitions $\{\lambda \mid \lambda_1
< n\}$ obtained by taking the {\it code} of the permutation;
see~\cite{Lam}.  If we make the identifications $\tF_w(x) =
\tF_{\lambda(w)}(x)$ and $s_w^{(k)}(x) = s_{\lambda(w)}^{(k)}(x)$
under this bijection, then this agrees with the usual indexing of
$k$-Schurs and affine Schurs by partitions; see~\cite{Lam}.
In~\cite{LLMS}, the monomial expansion of the $s_w^{(k)}(x)$ will be
given in terms of certain tableaux.

\subsection{Non-commutative $k$-Schur functions}
\label{sec:noncom}

The following definition was inspired by work of Fomin and
Greene~\cite{FG}.  Recall the isomorphism $\psi: \Lambda_n \simeq
\B$ of Theorem~\ref{thm:sym}.

\begin{definition}
Let $w \in W^0$.  The {\it non-commutative $k$-Schur functions} are
given by $$s^{(k)}_w:= \psi(s^{(k)}_w(x)) \in \B.$$
\end{definition}

\begin{example} \label{ex:kschur} Let $n = 3 = k+1$.  The affine Grassmannian
permutations with length less than or equal to 3 are $\id, s_0,
s_1s_0, s_2s_0, s_2s_1s_0$ and $s_1s_2s_0$.  Here we have computed
the corresponding non-commutative $k$-Schur functions, where we
write $A_{i_1i_2 \cdots i_k}$ for $A_{i_1}\cdots A_{i_k} =
A_{s_{i_1}s_{i_2}\cdots s_{i_k}}$.
\begin{align*}
s^{(k)}_{\id} &= 1 \\
s^{(k)}_{s_0} &= h_1 = A_0 + A_1 + A_2 \\
s^{(k)}_{s_1s_0} & = h_2 = A_{02}+A_{21} + A_{10} \\
s^{(k)}_{s_2s_0} & = h_1^2 - h_2 = A_{20}+A_{12}+A_{01} \\
s^{(k)}_{s_2s_1s_0} &= h_2h_1 = h_1h_2 = A_{021}+A_{010}+
A_{102}+A_{121}+A_{202}+A_{210} \\
s^{(k)}_{s_1s_2s_0} &= h_1^3 - h_2h_1 =
A_{120}+A_{010}+A_{201}+A_{121}+A_{202}+A_{012}.
\end{align*}
Observe that all the coefficients in the above formulae are
non-negative.  Also note that $s^{(k)}_{w}$ contains only one term
corresponding to a Grassmannian permutation and this term is $A_w$.
These properties are explained in Proposition~\ref{prop:kschur} and
Corollary~\ref{cor:stan}.
\end{example}

 The main result we need concerning the
non-commutative $k$-Schur functions is the following.  Define the
coefficients $b_{w,v} \in \Z$, where $w \in W^0$ and $v \in W_\aff$
by
\begin{equation}
\label{eq:kschur} s^{(k)}_w = \sum_{v \in W_\aff} b_{w,v} A_v.
\end{equation}

\begin{proposition}[\cite{Lam}, Proposition 42]
\label{prop:kschur} The coefficient $b_{w,v}$ is equal to the
coefficient of $\tF_{w}$ in the affine Stanley function $\tF_{v}$,
when written in the basis $\{ \tF_{w} : w \in W^0 \}$ of affine
Schur symmetric functions.  In particular, $s^{(k)}_w$ has a unique
Grassmannian term $A_w$.
\end{proposition}

In~\cite{Lam}, we conjectured that $b_{w,v}$ were positive and
showed how they implied positivity conjectures related to the toric
Schur functions of Postnikov~\cite{Pos} and $k$-Schur
functions~\cite{LLM}.  We shall return to these coefficients in
Section~\ref{sec:pos}.

\smallskip

We now describe the Hopf-structure on $\B$, acquired under the
isomorphism $\psi: \Lambda_n \simeq \B$.  We let the unit be $h_0 =
A_\id$ and the counit $\counit(b)$ to be the coefficient of $h_0$ in
$b$ when written as a polynomial in the $h_i$.

The coproduct $\Delta_\B: \B \rightarrow \B \otimes_\Z \B$ is given
by
$$
\Delta_\B(h_i) = \sum_{j \leq i} h_j \otimes h_{i-j}
$$
and extending $\Delta_\B$ to a ring homomorphism.

Now define $\omega_\B: \B \rightarrow \B$ by $\omega_\B(h_i) =
\psi(e_i(x))$. Let $\bar{}: \A_0 \rightarrow \A_0$ be the algebra
involution of the nilCoxeter algebra induced by $\bar{} : A_i
\mapsto A_{n-i}$.  The following result is
essentially~\cite[Proposition 16]{Lam}.

\begin{proposition}
Let $b \in \B$.  Then $\omega_B(b) = \bar b$.
\end{proposition}

For example, with $n = 3$ as in Example~\ref{ex:kschur},
$\omega_B(s^{(k)}_{s_1s_0}) = s^{(k)}_{s_2s_0}$.  We give $\A_0$ a
grading (half of the topologically induced one) by letting ${\rm
deg} A_i = 1$. Thus for an homogeneous element $b \in \B$ with
degree $d$ (as an element of $\A_0$) we define the antipode to be
$c_\B(b) = (-1)^d \; \bar b$, agreeing with the antipode of the
symmetric functions.

\section{Affine Grassmannian Schubert polynomials}
\label{sec:main} Our main result is the following.

\begin{theorem}\label{thm:main} The map $\theta: H_*(\G/\p) \to \Lambda_n$ given by
$$\theta: \sigma_w \longmapsto s_w^{(k)}(x)$$ is an isomorphism of
Hopf-algebras.  The map $\theta': H^*(\G/\p) \to \Lambda^n$ given by
$$\theta': \sigma^w \longmapsto \tF_w(x)$$ is an isomorphism of
Hopf-algebras.
\end{theorem}
In the homology case, this theorem was a conjecture of Mark
Shimozono, who based it on computer calculations using the nilHecke
ring and Theorem~\ref{thm:kk}.  The conjecture in the cohomology
case was made precise by Jennifer Morse.

Theorem~\ref{thm:main} will follow from two technical computations
comparing $\B$ and $\B'$ as Hopf algebras, which we now state.

\begin{proposition}
\label{prop:annihilate} Let $b \in \B$ and $s \in S$.  Then
$$\phi_0(bs) = \phi_0(s)b.$$
In other words, we have $\B \subset \B'$.
\end{proposition}
\begin{proposition}
\label{prop:coproduct} The two coproducts $\Delta$ and $\Delta_\B$
agree on $\B$ up to specialisation at 0:
\[
(\phi_0 \circ \Delta)(b) = \Delta_\B(b),
\]
for $b \in \B$.

\end{proposition}

\begin{theorem}
\label{thm:BB} The two algebras $\B$ and $\B'$ are identical as
subalgebras of $\A_0$.  Furthermore the two Hopf-structures agree
and we have for each $w \in W^0$,
$$
\phi_0(j(\sigma_{(w)})) = s_w^{(k)}.
$$
\end{theorem}
\begin{proof}
By Proposition~\ref{prop:kschur}, the element $s_w^{(k)} \in \B$ has
a unique Grassmannian term $A_w$.  By
Proposition~\ref{prop:annihilate}, it lies in $\B'$, so by
Proposition~\ref{prop:bprime}, we have $\phi_0(j(\sigma_{(w)})) =
s_w^{(k)}$ and the elements $s_w^{(k)}$ span $\B'$ so we conclude
that $\B = \B'$. For the Hopf algebra structure,
Proposition~\ref{prop:coproduct} shows that the coproduct agrees.
The agreement of the unit and counit is easy to see.

Finally, the antipode of $\B'$ (described after
Theorem~\ref{thm:pet} for any simple Lie type) is given by $t
\mapsto t^{-1}$ which is induced by the map $A_i \mapsto
-A_{\omega(i)}$ for each $i \in I$, where $\omega$ is the diagram
automorphism of the finite Dynkin diagram and $\omega(0) = 0$.  In
type $A_{n-1}$, this is given by $A_i \mapsto -A_{n-i}$. This agrees
with the antipode $c_\B$ of $\B$ described in
Section~\ref{sec:noncom}.
\end{proof}

\begin{proof}[Proof of Theorem~\ref{thm:main}]
For homology the theorem follows immediately from
Theorem~\ref{thm:BB} via the composition $$H_*(\Gr) \longrightarrow
\B \simeq \B' \stackrel{\psi^{-1}}{ \longrightarrow} \Lambda_n,$$
and the statement for cohomology follows by duality.
\end{proof}

\begin{remark}
Note that Proposition~\ref{prop:coproduct}, together with
Proposition~\ref{prop:kschur}, already imply that the coproduct of
$k$-Schur functions agrees with the product in cohomology.  Briefly,
the fact that each non-commutative $k$-Schur has a unique
Grassmannian ``leading term'' allows one to obtain a coproduct
expansion agreeing with Theorem~\ref{thm:kk} for Grassmannian
permutations. This is the argument given in~\cite{LamAS}.
\end{remark}

\begin{remark}
We find two symmetric function interpretations of
Theorem~\ref{thm:main} and Theorem~\ref{thm:BB} rather curious.

First, the Hopf algebra structure of the symmetric functions (which
has been studied from many perspectives) now attains a topological
origin.  In the classical identification of the Schur functions as
Schubert classes of the Grassmannian, only the product (not
coproduct) of the symmetric functions is used.  In particular it is
interesting that the Hall inner product is now given an
interpretation as the pairing between homology and cohomology.

Secondly, the commutativity of $\B = \B'$ links the symmetry of
affine Stanley symmetric functions with the commutativity of
$H_*(\Gr)$, which in turn is due to the fact that $\Gr$ is a double
loop space. Note that as shown in~\cite{Lam}, the usual Stanley
symmetric functions $F_w$ are special cases of affine Stanley
symmetric functions, so we have now a ``topological'' explanation of
the symmetry of Stanley symmetric functions.
\end{remark}

\subsection{Proof of Proposition~\ref{prop:annihilate}: $\B'$ contains $\B$}

We show that $\phi_0(h_i \cdot \alpha_j) = 0$ for each $i$ and the
result follows since $\{h_i\}_{i=0}^{n-1}$ generate $\B$ and
$\{\alpha_j\}_{j \in \Z/n\Z}$ generate the elements of $S$ with
constant term 0. Without loss of generality we will assume that
$j=1$. Let $I \subset \Z/n\Z$ be of size $i$. We calculate
$\phi_0(A_I \alpha_1)$ explicitly.  In the following table we let
$[2,r]$ be the largest interval of its form (possibly empty)
contained in $I$ which contains $2$.  It is possible that $[2,r]$
contains $0$ but it cannot contain 1 (since then it will have size
$n$).  Also the subset $I'$ never contains any of $0,1,2$. The sums
over $a$ are always over $a \in [2,r]$. The
$(\mathtt{A})$,$(\mathtt{B})$,$(\mathtt{C})$ are for marking the
terms only, for later use.
\begin{center}
\begin{tabular}{|c|c|c|}
 \hline & $I$ & $\phi_0(A_I \alpha_1)$ \\
\hline $0,1 \notin I$ & $I' \cup [2,r]$ &   $-\sum_a A_{I-\{a\}} (\mathtt{A})$ \\
$1 \in I$ and $0 \notin I$  & $I' \cup [2,r] \cup \{1\}$ & $2A_{I-\{1\}}(\mathtt{A}) +\sum_aA_{I-\{a\}}(\mathtt{C})$ \\
$0 \in I$ and $1 \notin I$ & $I' \cup [2,r] \cup \{0\}$ & $-A_{I-\{0\}}(\mathtt{A}) -\sum_aA_{I-\{a\}}(\mathtt{B})$\\
$0,1 \in I$ & $I' \cup [2,r] \cup \{0,1\}$ & $-A_{I-\{0\}}(\mathtt{C}) +A_{I-\{1\}}(\mathtt{B})$\\
\hline
\end{tabular}
\end{center}
For example \begin{align*} A_{[2,r]}A_1A_0 \alpha_1 \\
&=A_{[2,r]}A_1((\alpha_1 +\alpha_0)A_0 - 1) \\
&= -A_{[2,r]}A_1 + A_{[2,r]}(-\alpha_1A_1A_0 + 2 A_0 +
(\alpha_1+\alpha_0)A_1A_0 - A_0) \\
&= -A_{[2,r]-\{0\}} + A_{[2,r]-\{1\}} + \alpha_0A_{[2,r]}A_1A_0.
\end{align*}
The $A_t$ factors for $t \in I'$ always commute with the roots
$\alpha_j$ which appear in these calculations.  Note also that for
the exceptional set $I^* = \Z/n\Z-\{1\}$ which appears in the third
case above, the formula is to be interpreted as
$$\phi_0(A_{I^*}\alpha_1) = -2A_{I^* - \{0\}} -\sum_{a \in I^* \; ; \; a
\neq 0} A_{I^* - \{a\}}.$$

One observes that the terms marked $(\mathtt{A})$ or $(\mathtt{B})$
or $(\mathtt{C})$ when grouped together cancel out.  We have:
$(\mathtt{A})$ corresponds to subsets $J$ of size $i-1$ such that
$J$ contains neither $1$ nor $0$; and $(\mathtt{B})$ corresponds to
subsets $J$ of size $i-1$ such that $J$ contains $0$ but not $1$;
and $(\mathtt{C})$ corresponds to subsets $J$ of size $i-1$ such
that $J$ contains $1$ but not $0$. Every such subset in say case
$(\mathtt{A})$ will appear in all 3 case $(\mathtt{A})$ terms.
 No other subsets (those containing both $0$ and $1$) appear in the
sum $\sum_I A_I \alpha_1$.

For example, the element $A_J$ with $J = [2,4] \cup [5,7]$ will
appear in $\phi_0(A_I \alpha_1)$ for $I = [2,7]$ or $[1,4] \cup
[5,7]$ or $\{0\}\cup[2,4]\cup[5,7]$.  The multiplicities will be
$-1$, $2$, and $-1$ respectively, which cancel out.

\subsection{Proof of Proposition~\ref{prop:coproduct}: Comparison of coproduct} \label{sec:proof}
We proceed by computing directly $\phi_0 \circ \Delta$ on the
generators $h_i$ of $\B$.
\begin{lemma}
\label{lem:h} We have
\[
\phi_0(\Delta(h_i)) = \sum_{0 \leq j \leq i} h_{j} \otimes h_{i-j}.
\]
\end{lemma}
\begin{proof}
In the following computations the indices of roots and the elements
$A_i$ are to be taken modulo $n$.  Let $\beta_i = -\alpha_i$ be the
negative simple roots. We use $\Delta(A_i) = A_i \otimes 1 + 1
\otimes A_i + A_i \otimes \beta_i A_i$.

Let $i_1,i_2,\ldots,i_l \in \Z/n\Z$ be a cyclically decreasing
sequence.  For convenience we assume that if $i_k = i_j +1$ then $ k
= j-1$.  We have
\begin{eqnarray*} &\Delta(A_{i_1}A_{i_2}\cdots A_{i_l}) = \prod_j
\Delta(A_{i_j})\\& \!\!\!\!\!\!= \left(A_{i_1} \otimes 1 + 1 \otimes
A_{i_1} + A_{i_1} \otimes \beta_{i_1} A_{i_1}\right) \cdots
\left(A_{i_l} \otimes 1 + 1 \otimes A_{i_l} + A_{i_l} \otimes
\beta_{i_l} A_{i_l}\right)
\end{eqnarray*}
Let us expand the product, by picking one of the three terms in each
parentheses.  Strictly speaking we cannot multiply within $\A
\otimes_S \A$ since it is not a ring.  Instead we are calculating
the action of $\A$ on $\A \otimes_S \A$ via the coproduct: for
example $\Delta(A_iA_j) = \Delta(A_i)\cdot(\Delta(A_j)\cdot (1
\otimes 1)) \in \A \otimes_S \A$.  Formally, for our purposes the
calculation proceeds like multiplication.

Because of the cyclically decreasing assumption, the only times we
encounter a factor looking like $A_{i_a} \beta_{i_b}$ (where $a <
b$) we have either
\begin{equation}
\label{eq:commute} A_{i_a} \beta_{i_b} = \beta_{i_b} A_{i_a}
\end{equation}
or we will have $a = b-1$ and $i_{a+1} = i_{a} - 1$ and
\begin{equation}
\label{eq:notcommute} A_{i_a} \beta_{i_a-1} = (\beta_{i_a-1} +
\beta_{i_a})A_{i_a} + 1.
\end{equation}
If (\ref{eq:commute}) ever occurs, then $\beta_{i_b}$ commutes with
all $A_{i_c}$ where $c < b$ and we may ignore the term since
eventually we will apply $\phi_0$.  Similarly, if
(\ref{eq:notcommute}) occurs, the contribution of the term involving
$\beta_{i_a-1}$ is 0 after applying $\phi_0$.

Also we perform the calculation
\begin{equation}
\label{eq:power} A_{i+1} (\beta_i)^m = \beta_{i+1}^m A_{i+1} +
\beta_{i+1}^{m-1} + \mathrm{other\,terms,}
\end{equation}
where the other terms involve $\beta_i$ on the left somewhere (and
would be killed by $\phi_0$ later).

Let $B$ and $C$ be two subsets of $[0,n-1]$ with total size equal to
$i \leq n-1$.  We will first describe how to obtain the term $A_B
\otimes A_C$ (which occurs in $h_{|B|} \otimes h_{|C|}$) from
$\Delta(h_i)$.  Define a sequence of integers (``current degree'')
$(\cd(i) : j \in \Z/n\Z)$ by $\cd(j) = \max_{0 \leq t \leq n-1}\{|B
\cap [j-t,j]|+|C \cap[j-t,j]|-t-1\}$.  Since $|B|+|C| < n$ it is not
hard to verify that we can find $j$ so that $\cd(j) = 0$ and $i
\notin B \cup C$.

We may assume that $j = 0$.  Let $B = \{b_1 > \cdots > b_g > 0)$ and
$C = \{c_1 > \cdots > c_h > 0 \}$.  Define a sequence
$(t_1,t_2,\ldots,t_{n-1}) \in \{L,R,B,E\}^{n-1}$ as follows (where
$E =$ empty, $L =$ left, $R =$ right and $B =$ both):
\[
t_j = \begin{cases} E & \mbox{if $\cd(j) = 0$ and $E \notin B \cup
C$}\\ L & \mbox{if $\cd(j) = 0$ and $E \in B$ but $E \notin C$} \\
R & \mbox{if $E \notin B$ and ($\cd(j) > 0$ or $E \in C$)} \\
B & \mbox{otherwise.} \\
\end{cases}
\]
Now let $I = \{ j \in [1,n-1] \mid t_j \neq E\} \subset [1,n-1]$.
Then $A_B \otimes A_C$ is obtained from $\Delta(A_I)$ by picking the
term $A_{i_s} \otimes 1$ if $t_{i_s} = L$, the term $1 \otimes
A_{i_s}$ if $t_{i_s} = R$ and $A_{i_s} \otimes \beta_{i_s} A_{i_s}$
if $t_{i_s} = B$.

However, when the situation of (\ref{eq:power}) occurs, one has to
make a further choice between the two terms.  The sequence of
integers $(\cd(i))$ tells us the current degree (in the second
factor of the tensor product) in $S$ of the term that we want to
pick whenever we encounter the situation of (\ref{eq:power}).

For example if $\cd(t) = 3$ and $\cd(t+1) = 3$ then $t+1 \in B$ or
$t+1 \in C$.  In the first case we will have $(A_{t+1} \otimes 1)
\cdot (a \otimes \beta_i^3 b) $, for some $a$ and $b$ not involving
$S$, and there is no further choice. In the second case we get
$$(1 \otimes A_{t+1}) \cdot (a \otimes \beta_i^3 b) = a \otimes \beta_{i+1}^3\,A_{t+1}b +a \otimes \beta_{i+1}^2\,b,$$
modulo terms involving $\beta_i$ on the right.  One must make a
further choice between $\beta_{i+1}^3A_{t+1}$ and $\beta_{i+1}^2$.
We pick the first term since we want $t+1 \in C$ and this agrees
with the degree being $\cd(t+1) = 3$.

Thus every term of the form $A_B \otimes A_C$ appears in the
expansion of $\phi_0(\Delta(h_i))$.  Conversely, every ``current
degree'' sequence arising from choosing terms in the expansion of
$\Delta(A_I)$ as described above corresponds to some pair of subsets
$B$ and $C$ with total size $|B| + |C| = |I|$.
\end{proof}

\begin{proof}[Proof of Proposition~\ref{prop:coproduct}]
From Lemma~\ref{lem:h}, we have $\Delta_\B(h_i) =
\phi_0(\Delta(h_i))$. Now let $a \in \B$ and $b \in \B$ and suppose
we have shown that $\Delta_\B(a) = \phi_0(\Delta(a))$ and
$\Delta_\B(b) = \phi_0(\Delta(b))$.  Let $\Delta(a) = \sum_{w,v} A_w
\otimes a_{w,v} A_v$ and $\Delta(b) = \sum_{x,y} A_x \otimes b_{x,y}
A_y$, where $a_{w,v}, b_{x,y} \in S$.  Then
\begin{align*}
\phi_0(\Delta(ab))&= \phi_0(\Delta(a) \Delta(b)) \\
&=\phi_0(\sum_{w,v,x,y}A_wA_x \otimes a_{w,v} A_v b_{x,y} A_y) \\
&=\sum_{w,v,x,y}A_wA_x \otimes \phi_0(a_{w,v}) A_v \phi_0(b_{x,y})
A_y &\mbox{by Proposition~\ref{prop:annihilate}.} \\
&=\phi_0(\Delta(a))\phi_0(\Delta(b)) \\
&=\Delta_\B(a)\Delta_\B(b) \\
&=\Delta_\B(ab).
\end{align*}
Since the $h_i$ generate $\B$ this completes the proof.
\end{proof}

\section{Positivity}
\label{sec:pos} The connection between symmetric functions and the
geometry of the affine Grassmannian established in
Theorem~\ref{thm:main} allows us to resolve some positivity
properties of $k$-Schur functions and affine Schur functions.
\subsection{Positivity of cohomology product}
Kumar~\cite{Kum} and Graham~\cite{Gra} have shown that the structure
constants of $H^T(\Gr)$ in the Schubert basis are non-negative
polynomials.  Thus we have

\begin{corollary}
The product structure constants for the affine Schur functions
$\{\tF_w(x) : w \in W^0\}$ are non-negative integers, or
equivalently the coproduct structure constants for the $k$-Schur
functions $\{s_w^{(k)}(x) : w \in W^0\}$ are non-negative.
\end{corollary}

The statement for $k$-Schur functions was a conjecture
in~\cite[p.4]{LLM}.

\subsection{Positivity of homology product}
Peterson~\cite{Pet} showed that the structure constants of
$H_*(\Gr)$ are special cases of structure constants for the quantum
cohomology $QH^*(K/T)$ of the flag manifold.  Since these latter
structure constants are known to enumerate certain curves in $K/T$,
they are necessarily non-negative.  Thus we have

\begin{corollary}
\label{cor:hom} The coproduct structure constants for the affine
Schur functions $\{\tF_w(x) : w \in W^0\}$ are non-negative
integers, or equivalently the product structure constants for the
$k$-Schur functions $\{s_w^{(k)}(x) : w \in W^0\}$ are non-negative.
\end{corollary}

Evidence for Corollary~\ref{cor:hom} had recently been given
in~\cite{LM05}.  It was shown directly and combinatorially that the
structure constants for $k$-Schur functions included as a special
case the Gromov-Witten invariants of the Grassmannian.

\smallskip

A special case of the homology structure constants can be described
rather simply and explicitly.

\begin{corollary}
\label{cor:factorize1} Let $t_\lambda$ be a translation by an
antidominant co-root $\lambda \in Q^-$.  Then
$$
s^{(k)}_{t_\lambda}(x)s^{(k)}_{x}(x) = s^{(k)}_{xt_\lambda}(x).
$$
\end{corollary}
\begin{proof}
Combine Proposition~\ref{prop:specialclass} with
Theorem~\ref{thm:main}.
\end{proof}

If we translate the labeling of $k$-Schur functions from affine
Grassmannian permutations to partitions (see~\cite{Lam}) we obtain
the following result.

\begin{corollary}
\label{cor:factorize} Let $\mu$ be a partition with no part greater
than $k$ ($=n-1$), and $R$ be a rectangle of the form $R =
(a^{n-a})$. Then
$$
s^{(k)}_\rho(x) = s^{(k)}_{R}(x) s^{(k)}_\mu(x),
$$
where $\rho = R \cup \mu$.
\end{corollary}

A combinatorial proof of Corollary~\ref{cor:factorize} was given
in~\cite{LM04}. Peter Magyar communicates to us that he has a purely
geometric explanation of Corollary~\ref{cor:factorize1}.

In fact the homology structure constants are special cases of the
coefficients $b_{w,v}$ occurring in the expansion of an affine
Stanley symmetric $\tF_w(x)$ in terms of affine Schur functions, as
explained in~\cite{Lam}. By Proposition~\ref{prop:kschur}, the
$b_{w,v}$ also occur in expressing $s^{(k)}_w =
\phi_0(j(\sigma_{(w)}))$ in the $A_v$ basis.  Peterson~\cite{Pet}
also showed that these coefficients were positive, and thus this
resolves the combinatorial Conjecture 38 of~\cite{Lam}.

\begin{corollary}
\label{cor:stan} The affine Stanley symmetric functions $\tF_w(x)$
expand positively in terms of affine Schur functions.
\end{corollary}

As a special case we obtain the positive expansion of usual Stanley
symmetric functions in terms of Schur functions, first shown
in~\cite{EG,LS85}.  A combinatorial algorithm intended to prove
Corollary~\ref{cor:stan} was developed in~\cite{LamS}.

\begin{remark}
\label{rem:stan} There is a more direct geometric interpretation of
the coefficients $b_{w,v}$.  Consider the map $p: \Gr \rightarrow
\G/\BB$ given by the composition $\Omega K \hookrightarrow LK
\twoheadrightarrow LK/T$. Then we have $p^*(\sigma_\B^v) = \sum_{v
\in W^0} b_{w,v} \, \sigma^w$.  In other words, under the
isomorphism of Theorem~\ref{thm:main} the affine Stanley symmetric
functions $\tF_w(x)$ are the pullbacks $p^*(\sigma_\B^v)$.  This
also gives a direct geometric explanation of Stanley symmetric
functions.
\end{remark}

\subsection{Other positivity properties}
The $k$-Schur functions are known or conjectured to have other
positivity properties which may have interesting geometric
explanations.  As explained in Section~\ref{sec:kschur}, the
Macdonald polynomials are conjectured to be positive in the
$k$-Schur basis, and this suggests a precise connection between the
homology of $\Gr$ and Cherednik algebras.  The $k$-Schur functions
are also conjectured to be $k+1$-Schur positive, which suggests
positivity properties for Schubert classes under the embedding
$\Omega SU(n) \hookrightarrow \Omega SU(n+1)$.  Finally, the
original $t$-analogue $s^{(k)}(x;t)$ suggests that there may be an
additional $\C^*$-equivariance to be considered.

\section{Further directions}
\label{sec:final} There are many potential generalizations of this
work, as for the usual Schubert polynomials (\cite{LS}), one may
generalize in the $K$-theoretic, quantum, equivariant and other Lie
type directions. We note here that the $K$-theoretic version of the
nilHecke ring was developed in~\cite{KK2} and a $K$-theoretic
version of the affine Schur functions (called {\it affine Stable
Grothendieck polynomials}) was defined in~\cite{Lam}.

It would also be interesting to compare our work with the viewpoint
taken in geometric representation theory, for example via the
geometric Satake isomorphism.  In particular, Ginzburg~\cite{Gin}
gives a description of the cohomology of $\Gr$ as the symmetric
algebra $U(\check{\mathfrak{g}}^e)$ on the centralizer of a
principal nilpotent $e$ in the Langlands dual Lie algebra $
\check{\mathfrak{g}}$.  The homology Schubert classes corresponding
to $G(\O)$-orbits, that is those of the form $\sigma_{t_{\lambda}}$
where $\lambda \in Q^-$ is anti-dominant, are described explicitly
in~\cite{Gin} as matrix coefficients on $\check{\mathfrak{g}}^e$. It
would be interesting to compare this formula with the $k$-Schur
functions.

As evidence for a possible connection, we note that by
Proposition~\ref{prop:specialclass} we have $j(\sigma_{(t_\lambda)})
= \sum_{w \in W/W_\lambda} A_{t_{w \cdot \lambda}}$, and one can
compute directly that this element lies not only in the centralizer
$Z_{\A_\aff}(S)$ but also in the center $Z(\A_\aff)$ of $\A_\aff$.
It seems natural to ask whether the center $Z(\A_\aff)$ is generated
by these elements and the $W_\aff$-invariants $S^{W_\aff}$.  As the
affine nilHecke ring is in a loose sense a deformation of the affine
Hecke algebra, this compares well with the identification of the
spherical Hecke algebra as the center of the affine Hecke algebra;
see~\cite{Gai}.  We can show that the center of the affine
nilCoxeter algebra is indeed generated by the elements
$j(\sigma_{(t_\lambda)})$.

\begin{proposition}
The center of the affine nilCoxeter algebra $\A_0$ (for any Lie
type) is spanned by the elements
$$
\sum_{w \in W/W_\lambda} A_{t_{w\cdot \lambda}}
$$
where $W_\lambda \subset W$ is the stabiliser of $\lambda \in
Q^\vee$.
\end{proposition}

\begin{proof}
We use the length formula (\ref{eq:length}) repeatedly.

Let $a \in Z(\A_0)$ lie in the center.  We may assume that $a$ is
homogeneous.  Suppose first that the coefficient of $A_{wt_g}$ is
non-zero for some $w \neq {\rm id}$.  Pick a regular element $h \in
Q^\vee$ so that the product $w t_g t_h$ is length adding; that is
$\ell(wt_g t_h) = \ell(wt_g) + \ell(t_h)$ so that  $A_{wt_g}A_{t_h}
= A_{wt_{g+h}}$. But $A_{t_h}A_{wt_g}$ is either 0 or equal to $A_{w
t_{w.h + g}}$.  Since $w.h \neq h$, picking $h$ sufficiently large
we can ensure that the coefficient of $A_{wt_{g+h}}$ in $A_{t_h}
\cdot a$ is 0.  This shows that no term of the form $A_{wt_g}$ for
$w \neq \id$ is present.

Now suppose that $A_{t_h}$ occurs in $a$.  Suppose $r_i \cdot h \neq
h$.  Then exactly one of $r_it_h$ and $t_h r_i$ has length
$\ell(t_h) + 1$.  Comparing $A_i \cdot a$ with $a \cdot A_i$ we
deduce that the coefficient of $A_{t_h}$ is equal to the coefficient
of $A_{t_{r_i\cdot h}}$ in $a$.
\end{proof}

\bibliographystyle{amsalpha}

\end{document}